# Numerical strategy on the grid orientation effect in the simulation for two-phase flow in porous media by using the adaptive artificial viscosity method

Xiao-Hong Wang[a], Meng-Chen Yue[a], Zhi-Feng Liu[a,1*], Wei-Dong Cao[b], Yong Wang[c], Jun Hu[a], Chang-Hao Xiao[a], Yao-Yong Li[a]

[a] *Department of Thermal Science and Energy Engineering, University of Science and Technology of China, Hefei, Anhui 230027, P. R. China*

[b] *Information Management Center of Shengli Oilfield, Sinopec, Dongying, Shandong 257000, P. R. China*

[c] *Exploration and Development Research Institute, Shengli Oilfield Company, SINOPEC, Dongying, Shandong 257015, P. R. China*

## Abstract

It is a challenge to numerically solve nonlinear partial differential equations whose solution involves discontinuity. In the context of numerical simulators for multi-phase flow in porous media, there exists a long-standing issue known as Grid Orientation Effect (GOE), wherein different numerical solutions can be obtained when considering grids with different orientations under certain unfavorable conditions. Our perspective is that GOE arises due to numerical instability near displacement fronts, where spurious oscillations accompanied by sharp fronts, if not adequately suppressed, lead to GOE. To reduce or even eliminate GOE, we propose augmenting adaptive artificial viscosity when solving the saturation equation. It has been demonstrated that appropriate artificial viscosity can effectively reduce or even eliminate GOE. The proposed numerical method can be easily applied in practical engineering problems.

---

* Corresponding author. E-mail: lzf123@ustc.edu.cn (Zhi-Feng Liu).

Modeling multi-phase flows of immiscible fluids in porous media is important in various scientific and industrial fields. This long-lasting domain of interest has been enriched by many others, ranging from flow in packed bed reactors, remediation of dense non-aqueous phase liquids in contaminated soils, CO2 sequestration, hydrology, wastes storage in landfills or subsurface formations, flow in biological tissues, drying of porous materials, gas-water management in fuel cells, flow in filters and membranes, chemical engineering, nuclear safety and more. Numerical simulators of solving multi-phase flow in porous media are the products of evolution over several decades and have been widely used in many engineering fields. However conventional numerical algorithms suffer from the Grid Orientation Effect (GOE) [1-4], which produces different results when the orientation of the computational grid is changed even for simple geometries under certain unfavorable conditions. GOE significantly reduces the reliability of numerical simulation results, and produce serious errors in the prediction of the issues including oil and gas recovery, carbon neutrality, carbon capture, utilization and storage (CCUS) and *et al.* Though much effort has been directed toward the development of models without this effect, GOE remains one of the most challenging problems yet to be solved.

Solving nonlinear partial differential equations whose solution involves discontinuity is a formidable challenge. For instance, in the calculation of shock waves in aerodynamics, there is a long-standing problem known as the carbuncle phenomenon. In steady state blunt body calculations, Roe's scheme sometimes admits a spurious solution in which a protuberance grows ahead of the bow shock along the stagnation line [5,6]. This undesirable numerical phenomenon can also be observed in the numerical simulation of interstellar flows in astrophysics [7]. The carbuncle problem is considered as a typical example of numerical instability relevant to the discontinuity, and is still one of the greatest unresolved problem of classical numerical schemes in computational fluid mechanics [8].

It is well known that under certain conditions, the Buckley-Leverett equation, which is a typical hyperbolic equation describing the saturation filed, leads to shock waves [9]. In this letter, we focus on GOE in the numerical reservoir simulation. Our perspective is that the numerical instability near fronts cause GOE. Spurious numerical oscillations are always accompanied by sharp fronts and if not adequately suppressed, these oscillations can lead to GOE. To reduce, or even eliminate GOE, we propose augmenting adaptive artificial viscosity when solving the

Buckley-Leverett equation.

For the sake of simplicity and clarity, we consider incompressible two-phase flows in homogeneous and isotropic 2D media in the absence of capillary pressure. The conservative law can be expressed as:

$$\phi \frac{\partial S_\alpha}{\partial t} + \nabla \cdot \mathbf{V}_\alpha = 0, \tag{1}$$

where $\phi$ represents porosity and $\alpha = $ o or $\alpha = $ w denotes the oil and water phase, respectively. The phase velocity ($\mathbf{V}_o$ and $\mathbf{V}_w$) can be determined by the generalized Darcy's law:

$$\mathbf{V}_\alpha = - K \lambda_\alpha \nabla P, \tag{2}$$

where $K$ is the absolute permeability of the porous media. The phase mobility $\lambda_\alpha$ is defined as $K_{r\alpha}/\mu_\alpha$ where the relative permeability $K_{r\alpha}$ is a function of the phase saturation $S_\alpha$, and $\mu_\alpha$ is the phase viscosity. By introducing the total velocity $\mathbf{V} = \mathbf{V}_o + \mathbf{V}_w$, and the fractional flux function $F = \lambda_w/\lambda_T$ with the total mobility $\lambda_T = \lambda_o + \lambda_w$, Eq. (1) and Eq. (2) can be rewritten in a dimensionless form:

$$\begin{cases} \nabla \cdot \mathbf{V} = \nabla \cdot \left[ (K_{rw} + \frac{K_{ro}}{M}) \nabla P \right] = 0 \\ \frac{\partial S}{\partial t} + \mathbf{V} \cdot \nabla F = 0 \end{cases}, \tag{3}$$

with $M = \mu_o/\mu_w$, which denotes displaced-displacing phase viscosity ratio. Notice here in Eq. (3), all the variables are dimensionless.

The radial flow problem is suitable for performing numerical experiment to characterize the GOE. In this scenario, a single injection well is positioned at the center of a circular region. The reservoir is initially oil-saturated. At the circular outflow boundary, a constant pressure is imposed. We specially consider piston-type displacement, where the severity of GOE is more pronounced. Here, piston-type displacement refers to a distinct interface between the two fluids during the displacement process, and the distinct interface moves forward like a piston. Piston-type displacement occurs when the entropy condition in the Buckley-Leverett equation is satisfied, which can be expressed as $\frac{F(S_u)-F(S)}{S_u-S} \geq \frac{F(S_u)-F(S_d)}{S_u-S_d}$ for any value of $S$ between the upstream saturation $S_u$ and the downstream saturation $S_d$ [10,11]. Since $S_u = 1$ and $S_d = 0$ in this test, we set the fractional flux function $F(S) = S^2$, and thus the phase mobility $\lambda_w = \frac{F}{\mu_o(1-F)+\mu_w F}$, $\lambda_o = \frac{1-F}{\mu_o(1-F)+\mu_w F}$. It is evident that the entropy condition is met here, resulting in a preserved shock. The simulated domain is discretized using a regular Cartesian grid, and two classical numerical schemes - namely the five-point (5P) scheme and the nine-point (9P) scheme respectively, are employed to perform the numerical simulation.

The 9P scheme was derived from overlaying two 5P schemes associated with two square grids rotated relative to each other by $\pi/4$, which is considered capable of alleviating GOE [12,13]. Figure 1 illustrates the calculated saturation fields under different viscosity ratio $M$ and different gird size $\Delta x$ ($\Delta y = \Delta x$ in this study) for both 5P and 9P schemes. Notably, consistent results are obtained from both schemes when considering a viscosity ratio $M = 1$. As grid size decreases, simulation results converge towards a deterministic solution. However, when considering a viscosity ratio $M > 1$, the GOE occurs. The discrepancy between results obtained from the 5P and 9P schemes indicates GOE. Moreover, larger values of viscosity ratio $M$ results in a more obvious GOE. Although the 9P scheme alleviates GOE compared to its 5P counterpart, it fails to achieve convergence towards a deterministic solution as the grid size decreases.

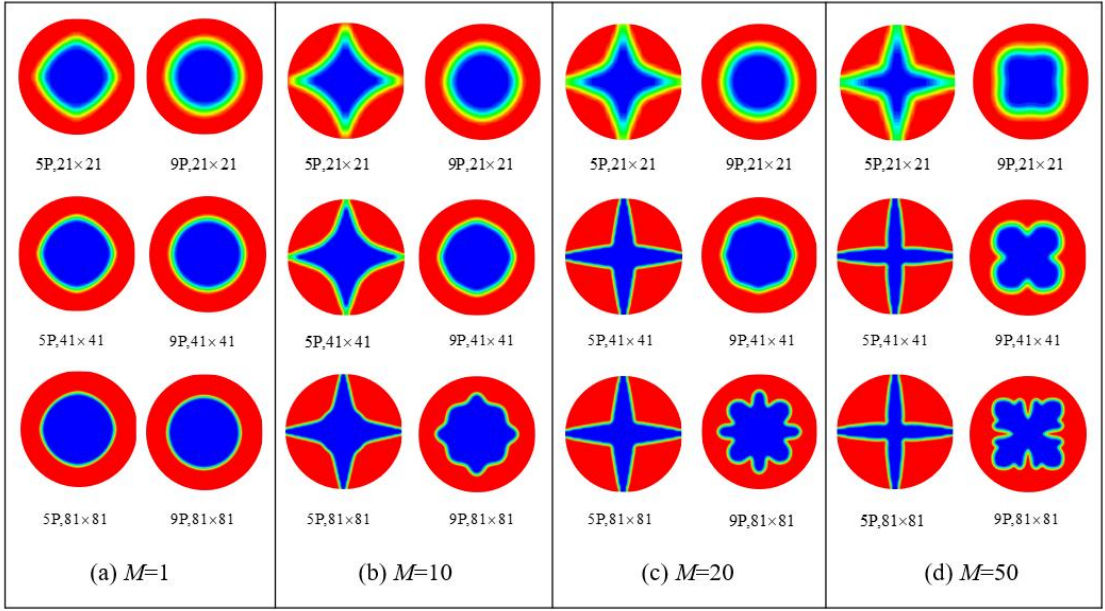

FIG. 1. The calculated saturation field of the radial flow from solving Eq.(3) (without artificial viscosity) under the different viscosity ratios $M$ and different numbers of grids ($21 \times 21$, $41 \times 41$ and $81 \times 81$, respectively) for the 5P and 9P schemes.

Starting from the perspective that GOE is a consequence of uncontrolled spurious numerical oscillations near the displacement front, we can introduce the artificial viscosity to suppress these oscillations and subsequently reduce or even eliminate GOE. The application of artificial viscosity method has proven effective in automatically eliminating wiggles behind shock fronts for hyperbolic equations. Artificial viscosity is a long-standing concept in computational fluid dynamics (CFD), successfully employed over past decades in numerous simulations involving fluid flows, such as solving hyperbolic problems like the Euler equations in gas dynamics [14-17]. The major difficulty in designing a highly accurate and robust artificial viscosity method is to make sure that a sufficient amount of stabilizing diffusion is added wherever it is needed, while in the rest of the computational domain the diffusion must be either switched off or small enough not to affect the high accuracy of the scheme there. At the same time, if the viscosity coefficient is too large in the areas away of discontinuity, the solution will be overly smeared there. Therefore, to achieve overall high resolution, the viscosity should be added in an adaptive way using a certain indicator, which should automatically pick rough parts of the computed solution and determine the amount of viscosity needed to be added there. In this letter, following Kurganov's idea [18], we augment the hyperbolic equation in Eq. (3) with an adaptive artificial viscosity, and then it becomes

$$\begin{cases} \nabla \cdot \mathbf{V} = \nabla \cdot \left[ (K_{\text{rw}} + \frac{K_{\text{ro}}}{M}) \nabla P \right] = 0 \\ \frac{\partial S}{\partial t} + \mathbf{V} \cdot \nabla F = C \nabla \cdot (\varepsilon \nabla S) \end{cases}, \quad (4)$$

where $C$ is a tunable positive viscosity coefficient and $\varepsilon = \varepsilon(S)$ is a nonnegative quantity whose size is automatically adjusted depending on the local properties of the discrete values of $S$. For computed solutions, $\varepsilon$ is set to be proportional to the size of the weak local residual (WLR), which serves as a smoothness indicator [19,20]. Denote the components of the total velocity $\mathbf{V}$ in the $x$ and $y$ direction as $u$ and $v$ respectively. The second equation in Eq. (4) can be rewritten as:

$$\frac{\partial S}{\partial t} + \frac{\partial f}{\partial x} + \frac{\partial g}{\partial y} = C \left[ \frac{\partial}{\partial x} \left( \varepsilon \frac{\partial S}{\partial x} \right) + \frac{\partial}{\partial y} \left( \varepsilon \frac{\partial S}{\partial y} \right) \right], \quad (5)$$

with $f = uF$ and $g = vF$.

Suppose a typical finite volume locates in $[x_{i-1/2}, x_{i+1/2}) \times [y_{j-1/2}, y_{j+1/2}) \times [t^{n-1/2}, t^{n+1/2})$, and the discrete variables are defined at the center points of the finite volume and denoted as $(x_i, y_j, t^n)$. Following Kurganov's algorithm, the 2D version of the WLR can be calculated as [18]:

$$E^{n-\frac{1}{2}}_{i+\frac{1}{2}, j+\frac{1}{2}} = \frac{1}{36\Delta} \Delta x \Delta y \mathcal{U}^{n-\frac{1}{2}}_{i+\frac{1}{2}, j+\frac{1}{2}} + \frac{1}{12\Delta} \left( \Delta y \Delta t \mathcal{F}^{n-\frac{1}{2}}_{i+\frac{1}{2}, j+\frac{1}{2}} + \Delta x \Delta t \mathcal{G}^{n-\frac{1}{2}}_{i+\frac{1}{2}, j+\frac{1}{2}} \right), \quad (6)$$

where

$$\Delta = \max(\Delta t, \Delta x, \Delta y), \quad (7)$$

$$u^{n-\frac{1}{2}}_{i+\frac{1}{2},j+\frac{1}{2}} = \left( S^n_{i+\frac{3}{2},j+\frac{3}{2}} - S^{n-1}_{i+\frac{3}{2},j+\frac{3}{2}} + S^n_{i+\frac{3}{2},j-\frac{1}{2}} - S^{n-1}_{i+\frac{3}{2},j-\frac{1}{2}} + S^n_{i-\frac{1}{2},j+\frac{3}{2}} - S^{n-1}_{i-\frac{1}{2},j+\frac{3}{2}} + S^n_{i-\frac{1}{2},j-\frac{1}{2}} \right.$$
$$\left. - S^{n-1}_{i-\frac{1}{2},j-\frac{1}{2}} \right) + 4\left( S^n_{i+\frac{3}{2},j+\frac{1}{2}} - S^{n-1}_{i+\frac{3}{2},j+\frac{1}{2}} + S^n_{i-\frac{1}{2},j+\frac{1}{2}} - S^{n-1}_{i-\frac{1}{2},j+\frac{1}{2}} + S^n_{i+\frac{1}{2},j+\frac{3}{2}} \right. \quad (8)$$
$$\left. - S^{n-1}_{i+\frac{1}{2},j+\frac{3}{2}} + S^n_{i+\frac{1}{2},j-\frac{1}{2}} - S^{n-1}_{i+\frac{1}{2},j-\frac{1}{2}} \right) + 16\left( S^n_{i+\frac{1}{2},j+\frac{1}{2}} - S^{n-1}_{i+\frac{1}{2},j+\frac{1}{2}} \right)$$

$$\mathcal{F}^{n-\frac{1}{2}}_{i+\frac{1}{2},j+\frac{1}{2}} = \left( f^n_{i+\frac{3}{2},j+\frac{3}{2}} - f^n_{i-\frac{1}{2},j+\frac{3}{2}} + f^n_{i+\frac{3}{2},j-\frac{1}{2}} - f^n_{i-\frac{1}{2},j-\frac{1}{2}} + f^{n-1}_{i+\frac{3}{2},j+\frac{3}{2}} - f^{n-1}_{i-\frac{1}{2},j+\frac{3}{2}} + f^{n-1}_{i+\frac{3}{2},j-\frac{1}{2}} \right.$$
$$\left. - f^{n-1}_{i-\frac{1}{2},j-\frac{1}{2}} \right) + 4\left( f^n_{i+\frac{3}{2},j+\frac{1}{2}} - f^n_{i-\frac{1}{2},j+\frac{1}{2}} + f^{n-1}_{i+\frac{3}{2},j+\frac{1}{2}} - f^{n-1}_{i-\frac{1}{2},j+\frac{1}{2}} \right) \quad (9)$$

$$\mathcal{G}^{n-\frac{1}{2}}_{i+\frac{1}{2},j+\frac{1}{2}} = \left( g^n_{i+\frac{3}{2},j+\frac{3}{2}} - g^n_{i+\frac{3}{2},j-\frac{1}{2}} + g^n_{i-\frac{1}{2},j+\frac{3}{2}} - g^n_{i-\frac{1}{2},j-\frac{1}{2}} + g^{n-1}_{i+\frac{3}{2},j+\frac{3}{2}} - g^{n-1}_{i+\frac{3}{2},j-\frac{1}{2}} + g^{n-1}_{i-\frac{1}{2},j+\frac{3}{2}} \right.$$
$$\left. - g^{n-1}_{i-\frac{1}{2},j-\frac{1}{2}} \right) + 4\left( g^n_{i+\frac{1}{2},j+\frac{3}{2}} - g^n_{i+\frac{1}{2},j-\frac{1}{2}} + g^{n-1}_{i+\frac{1}{2},j+\frac{3}{2}} - g^{n-1}_{i+\frac{1}{2},j-\frac{1}{2}} \right) \quad (10)$$

The fundamental concept of the adaptive artificial viscosity method is to set $\varepsilon$ as proportional to WLRs. In order to ensure that enough artificial numerical viscosity is augmented to rough parts of the solution, the discrete values of $\varepsilon$ can be selected as follows:

$$\varepsilon^n_{i+\frac{1}{2},j} = \max_{i,j}\left\{ \left|E^{n-\frac{1}{2}}_{i-\frac{1}{2},j-\frac{1}{2}}\right|, \left|E^{n-\frac{1}{2}}_{i+\frac{1}{2},j-\frac{1}{2}}\right|, \left|E^{n-\frac{1}{2}}_{i+\frac{3}{2},j-\frac{1}{2}}\right|, \left|E^{n-\frac{1}{2}}_{i-\frac{1}{2},j+\frac{1}{2}}\right|, \left|E^{n-\frac{1}{2}}_{i+\frac{1}{2},j+\frac{1}{2}}\right|, \left|E^{n-\frac{1}{2}}_{i+\frac{3}{2},j+\frac{1}{2}}\right| \right\}, \quad (11a)$$

$$\varepsilon^n_{i,j+\frac{1}{2}} = \max_{i,j}\left\{ \left|E^{n-\frac{1}{2}}_{i-\frac{1}{2},j-\frac{1}{2}}\right|, \left|E^{n-\frac{1}{2}}_{i-\frac{1}{2},j+\frac{1}{2}}\right|, \left|E^{n-\frac{1}{2}}_{i-\frac{1}{2},j+\frac{3}{2}}\right|, \left|E^{n-\frac{1}{2}}_{i+\frac{1}{2},j-\frac{1}{2}}\right|, \left|E^{n-\frac{1}{2}}_{i+\frac{1}{2},j+\frac{1}{2}}\right|, \left|E^{n-\frac{1}{2}}_{i+\frac{1}{2},j+\frac{3}{2}}\right| \right\}, \quad (11b)$$

Finally, the parameter $C$ in Eq. (5) is a positive viscosity coefficient that must be carefully chosen to effectively control the quality of the computed solution. The value of the parameter $C$ can be chosen as:

$$C = \frac{(\Delta x)^2 + (\Delta y)^2}{\alpha \Delta t \varepsilon^n_{\max}}, \quad (12)$$

with $\varepsilon^n_{\max} = \max_{i,j}\{\varepsilon^n_{i+1/2,j}, \varepsilon^n_{i,j+1/2}\}$ and a conservative choice being $\alpha = 4$. Here, the term 'conservative choice' implies that if $\alpha = 4$, an adequate amount of artificial viscosity will be augmented and correspondingly the fronts will be smeared.

This adaptive artificial viscosity implementation involves using diffusion coefficients that are very small in smooth regions, and becomes strong enough near piston-like fronts to prevent oscillations. Figure 2 shows the saturation fields obtained by solving Eq. (5) under different viscosity ratio $M$ and different gird size $\Delta x$ ($\Delta y = \Delta x$ in this study). In these calculations, the time step is set as $\Delta t = 33.5 \times \Delta x^{3.3}$. In order to obtain the fronts as narrow as possible, we set $\alpha = 67$ for 5P scheme and $\alpha = 400$ for 9P scheme in this test. It is amazing that for both 5P and 9P scheme, the irregular front shape disappears and is replaced by a circle which is consistent with the analytical solution. Here, the analytical solution can be expressed as $r_f = \sqrt{\frac{Qt}{\pi\phi}}$ with $r_f$ indicating the location of the circular front and $Q$ being the injection flow rate. Although the artificial viscosity widens and smear out the front, decreasing grid size leads to a convergent and sharp front.

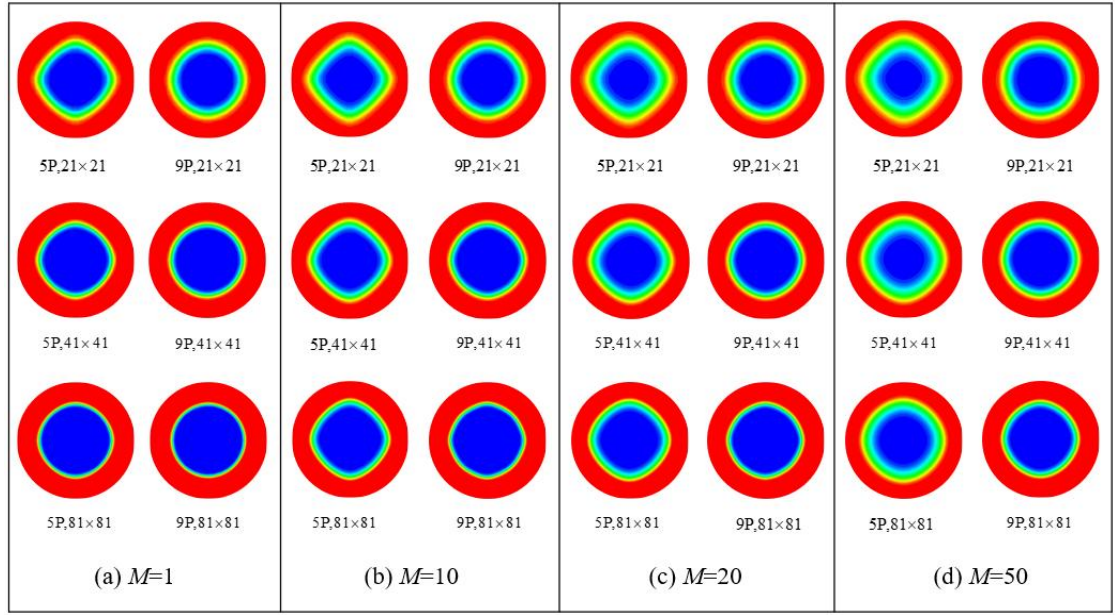

FIG. 2. The calculated saturation field of the radial flow from solving the Eq. (4) (with adaptive artificial viscosity) under different viscosity ratio $M$ and different numbers of grids ($21 \times 21$, $41 \times 41$ and $81 \times 81$, respectively) for the 5P and 9P schemes.

Another test is performed, considering a typical scenario with one injection well and four production wells (1IW-4PW). The time step is also set as $\Delta t = 33.5 \times \Delta x^{3.3}$, and the parameter $\alpha$ is assigned as $\alpha = 7$ for 5P scheme and $\alpha = 40$ for 9P scheme in this particular test. As shown in Fig. 3, when solving the original conservative Eq. (3) without artificial viscosity, significant GOE can be observed for high values of viscosity ratio $M$. However, by solving the conservative Eq. (4) with augmented artificial viscosity, the GOE is nearly eliminated, resulting in consistent outcomes between the 5P and 9P schemes (refer to Fig. 4).

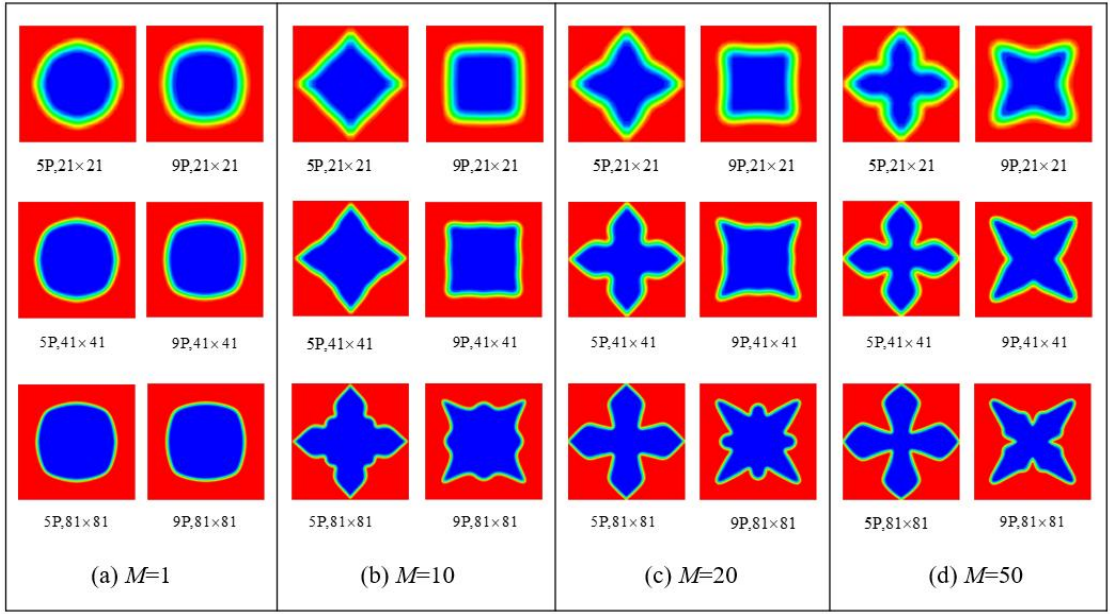

FIG. 3. The calculated saturation field of 1IW-4PW case from solving Eq.(3) (without artificial viscosity) under different viscosity ratio $M$ and different numbers of grids ($21 \times 21$, $41 \times 41$ and $81 \times 81$, respectively) for the 5P and 9P schemes.

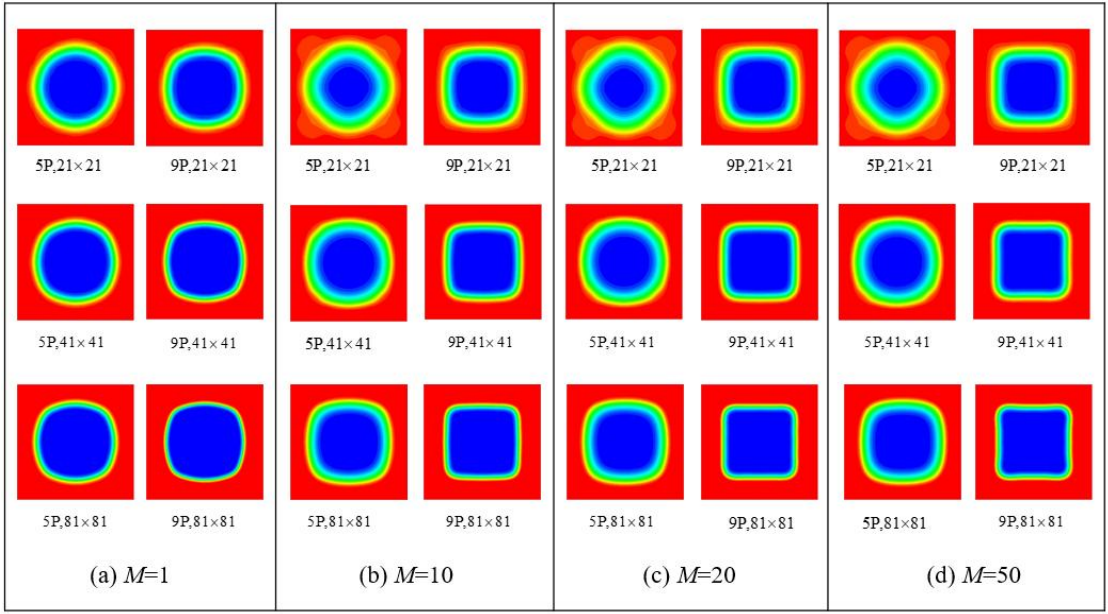

FIG. 4. The calculated saturation field of 1IW-4PW case from solving Eq. (4) (with adaptive artificial viscosity) under different viscosity ratio $M$ and different numbers of grids ($21 \times 21$, $41 \times 41$ and $81 \times 81$, respectively) for the 5P and 9P schemes.

In recent decades, extensive research has been carried out on the phenomenon GOE, including the utilization of irregular grids to mitigate its detrimental impact on practical engineering problems [20-24]. It has been demonstrated that the instability of 5P numerical scheme leads to GOE [3], but its underlying physical mechanism remains an enigmatic issue. While a linear stability analysis has been performed for multiphase flow in porous media, it was mainly focused on the convection-diffusion equation [25-28]. Chorin proposed that perturbations near the fronts contribute to increased instability; however, this stability analysis did not consider boundary condition constraints [29]. To the best of our knowledge, a rigorous instability analysis specifically addressing discontinuous solutions of pure nonlinear hyperbolic equations is still lacking. During carbuncle phenomenon research, it has been suggested that numerical schemes with higher artificial viscosity are advantageous in reducing this adverse occurrence. In our study on multiphase flow in porous media, we explicitly augment adaptive artificial viscosity to effectively reduce or even eliminate GOE. To enhance resolution when describing fronts, grid sizes need to be sufficiently small. As grid size approaches zero according to Eq. (12), the augmented artificial viscosity near the fronts tends towards infinity, which is different with the treatment of carbuncle phenomenon.

For dealing with the long-standing problem GOE, we apply the adaptive artificial viscosity method to solve two-phase flows of immiscible fluids in porous media. It is demonstrated that GOE can be effectively reduced, or even eliminated with properly added artificial viscosity. This approach is expected to be widely applied in various engineering fields such as oil and gas recovery, chemical engineering, nuclear safety, carbon neutrality, CCUS and *et al.*